\documentclass[a4paper,reqno,hidelinks]{amsart}

\pdfoutput=1

\usepackage{tikz}
\usetikzlibrary{calc}
\usetikzlibrary{arrows}
\usetikzlibrary{decorations.markings}
\usetikzlibrary{shapes} 
\usetikzlibrary{fit}
\usepackage{amsmath}
\usepackage{amsfonts}
\usepackage{amssymb}
\usepackage{hyperref}
\usepackage{cleveref}
\usepackage{mathtools}
\usepackage{enumitem}

\usepackage[llparenthesis,rrparenthesis]{stmaryrd}

\newcommand{\type}{\mathfrak C}
\newcommand{\UU}{\ensuremath{\mathcal{U}}}

\newcommand{\ty}{\mathsf{T}}
\newcommand{\fst}{\mathsf{shape}}
\newcommand{\unit}{\mathbf{1}}
\newcommand{\fullnerv}{N}

\newcommand{\plus}[1]{{#1}_{\scalebox{0.7}{+}}}
\newcommand{\minus}[1]{{#1}_{\scalebox{0.95}{-}}}
\newcommand{\posnerv}{\plus N}
\newcommand{\deltplus}{\plus{\Delta}}

\newcommand{\oppo}[1]{{#1}^{\scalebox{0.7}{\text{op}}}}

\newcommand{\opandplus}[1]{{#1}_{\scalebox{0.7}{+}}^{\scalebox{0.7}{\text{op}}}}
\newcommand{\deltop}{\opandplus{\Delta}}

\usepackage{mathtools}
\newcommand{\jdeq}{\equiv}
\newcommand{\defeq}{\vcentcolon\equiv}

\newcommand{\fib}{\twoheadrightarrow}

\newcommand{\mapfunc}{\mathsf{ap}}

\newcommand{\ct}{%
  \mathchoice{\mathbin{\raisebox{0.5ex}{$\displaystyle\centerdot$}}}%
             {\mathbin{\raisebox{0.5ex}{$\centerdot$}}}%
             {\mathbin{\raisebox{0.25ex}{$\scriptstyle\,\centerdot\,$}}}%
             {\mathbin{\raisebox{0.1ex}{$\scriptscriptstyle\,\centerdot\,$}}}
}

\newcommand{\Nop}{\mathbb{N}^{\mathsf{op}}}
\newcommand{\N}{\mathbb{N}}
\newcommand{\E}{\mathbf{E}}

\newcommand{\uu}{{\mathrm u}}
\newcommand{\vv}{{\mathrm v}}
\newcommand{\ww}{{\mathrm w}}
\newcommand{\sss}{{\mathrm s}}
\newcommand{\pp}{{\mathrm p}}
\newcommand{\qq}{{\mathrm q}}

\newcommand{\id}{\mathsf{id}}
\newcommand{\co}{\colon}

\newcommand{\term}{\mathsf{1}}

\newcommand{\istype}{\;\, \textit{type}}
\newcommand{\Sp}{\mathsf{Sp}}
\newcommand{\refl}{\mathsf{refl}}
\newcommand{\D}{\mathfrak D}
\renewcommand{\deg}{\mathsf{deg}}

\renewcommand{\fib}{\twoheadrightarrow}
\newcommand{\lfib}{\twoheadleftarrow}
\newcommand{\cof}{\rightarrowtail}

\def\defthm#1#2#3{
  \newtheorem{#1}[theorem]{#2}
  \newtheorem*{#1*}{#2}
  \crefname{#1}{#2}{#3}
}

\newtheorem{theorem}{Theorem}[section]
\newtheorem*{theorem*}{Theorem}
\crefname{theorem}{Theorem}{Theorems}

\defthm{proposition}{Prop.}{Props.}
\defthm{lemma}{Lemma}{Lemmata}
\defthm{corollary}{Cor.}{Cors.}

\theoremstyle{definition}
\defthm{remark}{Rem}{Rems.}
\defthm{definition}{Def.}{Defs.}
\defthm{example}{Ex.}{Exs.}

\crefname{section}{Sec.}{Secs.}
\crefname{subsection}{Subsec.}{Subsecs.}
\crefname{figure}{Fig.}{Figs.}

\DeclarePairedDelimiter\braces\lbrace\rbrace
\DeclarePairedDelimiterX\set[2]\lbrace\rbrace{#1 \mathrel{}\mathclose{}\delimsize|\mathopen{}\mathrel{} #2}

\title{Space-Valued Diagrams, Type-Theoretically
\\
(Extended Abstract)
}

\author{Nicolai Kraus \and Christian Sattler}
\thanks{Grant acknowledgments: Nicolai Kraus is supported by EPSRC grant EP/M016994/1, and Christian Sattler by the Air
Force Research Laboratory, under agreement number FA8655-13-1-3038.}

\begin{document}

\begin{abstract}
Topologists are sometimes interested in space-valued diagrams over a given index category,
but it is tricky to say what such a diagram even is if we look for a notion that is stable under equivalence.
The same happens in (homotopy) type theory, where it is known only for special cases how one can define a type of type-valued diagrams over a given index category.

We offer several constructions.
We first show how to define homotopy coherent diagrams which come with \emph{all} higher coherence laws explicitly, with two variants that come with assumption on the index category or on the type theory. 
Further, we present a construction of diagrams over certain Reedy categories. 
As an application, we add the degeneracies to the well-known construction of semisimplicial types, yielding a construction of \emph{simplicial types} up to any given finite level.

The current paper is only an extended abstract, and a full version is to follow.
In the full paper, we will show that the different notions of diagrams are equivalent to each other and to the known notion of Reedy fibrant diagrams whenever the statement makes sense.
In the current paper, we only sketch some core ideas of the proofs.
\end{abstract}

\maketitle

\section{Introduction and Background} \label{sec:motivation}

Working with categorical constructions internally in (homotopy) type theory is a delicate task. 
To illustrate this, consider the possibly most straightforward approach of defining the notion of a category, which might be the following.
We may say that a category consists of a type $A$ of objects; for any $x,y : A$, a type of morphisms $\mathsf{Hom}(x,y)$; a composition operation $\circ : \mathsf{Hom}(y,z) \times \mathsf{Hom}(x,y) \to \mathsf{Hom}(x,z)$; and the identities.
On top of this, we would like to add some laws, namely associativity $h \circ (g \circ f) = (h \circ g) \circ f$ and rules for the identities.
However, when stating these rules in type theory, our only option is to use the internal equality type.
This does not simply give us proof-irrelevant laws, but it gives us new structure.\footnote{%
Remark: In some formulations of dependent type theory, the axiom \emph{UIP}, \emph{uniqueness of identity proofs}, sometimes called Streicher's \emph{axiom K}, is assumed.
This ensures that any two elements of the same equality type are themselves equal. 
In this case, one \emph{can} treat equalities as laws rather than structure. 
This axiom changes many aspects of the theory, and many parts of our current paper would (just as many other results of homotopy type theory) significantly simplify in such a setting.
In the current paper, we consider type theories without \emph{UIP}, in particular theories considered in homotopy type theory, where \emph{UIP} would be an inconsistent axiom. (Our strictification construction could be of interest even in a system with \emph{UIP}.)%
}
To take an example, there is more than one way to prove an equality of the form $k \circ (h \circ (g \circ f)) = ((k \circ h) \circ g) \circ f$.
Put differently, the associativity rule cannot be treated as a law, but rather has to be seen as an operator.
To ensure this operator's well-behavedness, one needs to add a rule which corresponds to the well-known pentagon one has in the definition of a bicategory.
Unfortunately, this pentagon can again not be stated as a proof-irrelevant law, but constitutes structure requiring coherences.
If we continue this way, what we might eventually get is an $(\infty,1)$-category, but stating the full definition in this way is an open problem.

As stating all coherences is so involved,
one may suggest to simply ignore the need for coherence structure and laws at some point and settle for an ``incomplete'' definition. 
The downside of such an approach is that certain constructions will not work or lead to structure weaker than the one we have started with.
For example, consider an (ordinary) category $C$ with an object $x$, for which we have the slice category $C \slash x$ (sometimes also written as comma category $C \downarrow x$). 
One can quickly check that composition in the slice category requires the associativity law in $C$.
One level higher, in the case that $C$ is a bicategory, we do need coherence of associativity in $C$ (i.e.\ the pentagon law) in order to show that composition in $C\slash x$ is associative, and this ``level-shift'' seems to happen at higher levels in the same way.

What is particularly intriguing is that type theory \emph{does} have some structures which form an ordinary ``strict'' category, where the associativity and identity laws hold \emph{judgmentally} (i.e.\ both sides of the relevant equation are definitionally equal, i.e.\ have the same normal forms).
The standard example is a type universe: objects of the corresponding category are types, and morphisms are functions. 
Associativity of function composition and identity laws hold ``on the nose'' (at least under the usual assumption of a judgmental $\eta$-law for functions).
However, being unable to state what such a strict category is inside type theory, we are of course unable to describe this property internally.

Without the possibility to express the strict categorical structure of a universe $\UU$ in type theory,
it is also difficult to say what a type-valued diagram over a given index category $\mathcal D$ is.
Say, $\mathcal D$ is an externally fixed category (i.e.\ it is given as a normal category outside of the type theory), and we want to express inside type theory what a type-valued diagram over $\mathcal D$, i.e.\ a  functor from $\mathcal D$ to $\UU$, is.
In other words, we want to find a type of such functors.
Again, a somewhat canonical attempt is to say that, to give a functor $F$, we should have a type $F(X) : \UU$ for every object $X$ in $\mathcal D$; a function $F(f) : F(X) \to F(Y)$ for every morphism $f \in \mathcal D(X,Y)$; and equalities for the functor laws.
Unsurprisingly, we get the same problem as before: the equalities we give form structure which requires coherence in order to be well-behaved, but the tower of needed coherence laws is a priori infinite and hard to describe. 

At this point, we want to note that the described phenomena are not at all new discoveries of type theory.
Quite the contrary, these problems have been well-known in the mathematical communities for algebraic topology / homotopy theory for a rather long time. 
There, people have encountered very similar problems, namely that certain properties should be expressed in a homotopy-invariant way in order to be well-behaved, which however requires coherence conditions that are hard to handle (see Boardman--Vogt~\cite{bv:htpy-invar}).
What (homotopy) type theory does is offering a new view on the same old problems.
By offering a new view, it also has the potential to offer new approaches which can be used in the original mathematical settings. 
In our case, ``type-valued diagrams'' could for example be translated quite naively to ``space-valued diagrams'', for a notion of space that is modelled by something close to a fibration category as introduced by~\cite{brown:fibration-categories}.
In this context, the work of \cite{cisinski:fibration-categories,rb:fibration-categories,szumilo2014two,KapulkinSzumilo2016} is very related.
For example, \cite[Lem.~3.18]{szumilo2014two} gives an alternative correspondence between what we call weak diagrams with and without identities.

The two described problems, namely describing the structure of a category, and describing type-valued functors,
are very closely related.
If we were able to express the coherent categorical structure of $\UU$, we could reasonably hope to be able to use it for expressing the coherent structure of functors into $\UU$.
Vice versa, one model for $(\infty,1)$-categories (that is, categories with infinitely many ``levels of coherence structures'') are \emph{complete Segal spaces}~\cite{rezk2001model};
and indeed, these are space-valued diagrams with certain conditions over the category $\Delta^\mathsf{op}$.
Preliminary suggestions have been made to translate these concepts to type theory~\cite{altenCapKraus_infCats}, and this line of research is currently work in progress 

For some specific (externally given) index categories $\mathcal D$, it is known how to represent type-valued diagrams over $\mathcal D$.
Let us give an overview:

\begin{enumerate}[leftmargin=*]
\item {\textbf{Categories generated by simpler structure:}} \label{item:nolaws}
 If $\mathcal D$ is a finite discrete category, then the corresponding diagrams can be described as finite produces $\UU \times \ldots \times \UU$. 
 More generally, $\mathcal D$ could be generated by something that has objects and morphisms, but fewer laws, for example a finite directed graph.
 
\item {\textbf{Groupoidal index categories:}} \label{item:groupoidal}
 The following special case was pointed out to us by Steve Awodey, Ulrik Buchholtz, and Peter LeFanu Lumsdaine. 
 If $\mathcal D$ happens to be a finite groupoid, where all morphisms are invertible, 
 one can form a type $K(\mathcal D)$ representing $\mathcal D$, along the lines of Eilenberg-MacLane spaces~\cite{licataFinster_Eilenberg}, using higher inductive types.
 The type $K(\mathcal D)$ can be seen as an internal representation of $\mathcal D$, and diagrams over $\mathcal D$ will then simply be given by the function type
 $K(\mathcal D) \to \UU$. 
 
 If we were working in some form of \emph{directed homotopy type theory}, another topic of current research, 
 we speculate that we could do the same for a not necessarily groupoidal category $\mathcal D$, by using a directed version of higher inductive types.

\item {\textbf{Internally represented categories:}} 
 In the two cases \eqref{item:nolaws} and \eqref{item:groupoidal}, 
 the category $\mathcal D$ can sometimes (possibly partially) be represented internally, and when this happens, it usually helps us to relax the ``finiteness'' condition.
 For example, when the objects of a discrete category are given as a not necessarily finite type (set) $A$, we can just take the function type $A \to U$. 
 For graphs, this is discussed in a textbook exercise~\cite[Ex.~7.2]{HoTTbook}.

\item {\textbf{Truncated types:}}
 Instead of simplifying the index category, we can simplify the involved types. 
 More precisely, we can consider the situation in which all occurring relevant types are truncated at a specific level.
 The precise technical meaning of this is not important for our current paper, but the interested reader is invited to consult the homotopy type theory textbook~\cite[Chap.~7]{HoTTbook}.
 This is a slight generalisation of the situation that \emph{UIP} is assumed, which would represent the $0$-truncated case.
 It is essentially the approach chosen in the development of \emph{univalent category theory} by Ahrens, Kapulkin, and Shulman~\cite{ahrens_Rezk}: the type of objects of a category in their sense is required to be $1$-truncated.
 This allows them to cut off the required towers of coherence laws at a very low level and develop concepts from standard category theory neatly.
 Note however that the truncation condition means that the universe $\UU$ is not a category in their sense, meaning that this attempt does not work for us.

\item {\textbf{Strict equality:}}
 We could consider a theory with a notion of strict equality in the style of Voevodsky's \emph{homotopy type system}, called \emph{HTS}~\cite{voe_HTS}, or another form of two-level theory~\cite{altCapKra_twolevels}.
 With strict equality at hand, the categorical and functor laws could be formulated directly, without the requirement of coherence laws.
 This however is not what we want here. 
 The resulting constructions would not be homotopy invariant, i.e.\ would not be \emph{fibrant types}, and all the nice properties of types in homotopy type theory would be lost.
 This is not to say that a two-level system is useless here: it can serve as a tool to formulate and streamline, and maybe even implement, certain concepts, as will be described in forthcoming work by Annenkov, Capriotti and a current author~\cite{Annenkov_two}.
 However, it does not replace any of the work which needs to be done here. 
 
\item {\textbf{The inductive approach (Reedy fibrant diagrams over inverse categories):}}
 A further special case is the one where $\mathcal D$ is an inverse category, meaning that there is some form of well-founded ordering on the objects, and non-identity morphisms only go from larger objects to smaller objects. 
 In this case, we prefer to write $I$ instead of $\mathcal D$.
 The point of an inverse category $I$ is that a type-valued functor on $I$ can be encoded in an inductive fashion, using the given ordering.
 This way of representing (strict) diagrams over an inverse category has been explored in detail by Shulman~\cite{shulman_inversediagrams}
 and is of great importance for our current paper. 
 Thus, we will review the constructions in detail in \cref{sec:introtoreedy} below.
\end{enumerate}
 
 If we are given $I$, we can \emph{externally} consider Reedy fibrant diagrams over $I$ (which are simply strict functors from $I$ into a category of types, with a certain property).
 We can also try to \emph{internalise} everything and define a \emph{type} of such diagrams inside type theory (one benefit of which is that it allows to implement the construction in a proof assistant).
 Whether this internalisation of diagrams over $I$ is possible in full depends on both $I$ and the type theory.
 If $I$ is an (externally fixed) \emph{finite} inverse category, the type of Reedy fibrant diagrams over $I$ can always be encoded in ``standard'' homotopy type theory, see Shulman's work~\cite{shulman_inversediagrams}.
 If $I$ is infinite, this is in general not the case. 
 (Note that by \emph{``standard'' homotopy type theory} we mean the type theory presented in the textbook~\cite{HoTTbook}.)
 
 For the very concrete example that $I$ is the opposite of the category of finite no-nempty ordinals and increasing injective functions, written $\deltop$, 
 the corresponding Reedy-fibrant diagrams are well-known as \emph{semisimplicial types} in the community of homotopy type theory.
 Whether this construction can be internalised, i.e.\ whether we can write down a type of semisimplicial types, is a well-known open problem.
 It is known that we can do it if we restrict ourselves to $\deltop$ up to a fixed finite level (as this would be diagrams over a finite inverse category).
 For the unrestricted (or restricted only by an internal variable) category $\deltop$, the answer to the problem is unknown.
 It has been subject of numerous informal discussions and is, 
 for example, recorded in~\cite{herbelin_semisimpl,nicolai:thesis,shulman:eating}.
 Our inability to perform such a construction is by some people seen as a major incompleteness of homotopy type theory as presented in~\cite{HoTTbook}.
 This has triggered the development of more powerful type theories with strict equalities as discussed above, where the encoding is possible.
 The original suggestion in this direction is Voevodsky's \emph{HTS}~\cite{voe_HTS}, one version of which has been made precise in~\cite{altCapKra_twolevels}.
 What exactly is needed to make the internal construction of semisimplicial types possible is content of active research, but a popular and sufficient assumption is that the natural numbers of the ``strict'' fragment of the theory coincide with the ordinary (``fibrant'') type natural numbers.
 Another (similar) suggestion of a theory that allows the construction of semisimplicial type was made by Part and Luo~\cite{DBLP:journals/corr/PartL15}.

 A further problem sometimes discussed is whether and how one can add \emph{degeneracies} to obtain \emph{simplicial types}; this is hard even if we restrict ourselves to a finite part of $\deltop$, as the index category is in this case not inverse anymore. 
 Semisimplicial types are important for the current paper and will be discussed in more detail in \cref{sec:introtoreedy}.

\subsection*{Contributions}  
In this paper, we present the following.

\paragraph*{\textbf{Homotopy coherent} (or simply \textbf{weak}) \textbf{diagrams}}
 For a finite inverse category $I$, we make precise the idea of defining a type-valued diagram as the collection of a type for each object, 
 a function for each morphism, a composition operator for each composable pair of morphisms, and \emph{all} the coherence conditions.
 As we will see, the (a priori infinite) tower of coherence laws can be expressed in a finite way and can thus be internalised in ``standard'' homotopy type theory (i.e.\ it can be implemented in a proof assistant). 
 This finiteness property holds thanks to two conditions: first, the inverse property of $I$ allows us to essentially remove identities, viewing $I$ as freely generated from an inverse semicategory, and second, $I$ being finite and inverse means that the coherence laws for ``different orders of compositions of arrows'' are only needed up to a certain level.
 We call these diagrams \emph{homotopy coherent} or, for the sake of convenience, simply \emph{weak}, while we refer to Reedy fibrant diagrams as \emph{strict}.
 
 Unsurprisingly, a strict diagram gives rise to a weak one almost directly. 
 We present a \emph{strictification construction} for the other direction, which turns a weak diagram into a strict Reedy fibrant diagram.
 
 If we internalise and look at the \emph{types} of weak and strict diagrams, our strictification establishes an equivalence.
 The proof of this statement is only sketched in a later part (\cref{subsec:weakStrictEquiv}) of the current paper.
 It requires technical constructions and machinery that will be presented in a future full version of this paper~\cite{we_futurework}. 
  
\paragraph*{\textbf{General homotopy coherent diagrams}}
 In the next step, we generalise our construction of homotopy coherent diagrams over inverse categories to an arbitrary index category $C$.
 To do this, we put back identities.
 We do this in an economic way which follows a general principle suggested by Paolo Capriotti, who in turn was inspired by Harpaz~\cite{harpaz2015quasi}.
 This gives us a notion of general homotopy coherent diagrams over any category $C$.
 
 The downside is that such a diagram will nearly always consist of infinitely many components.
 It is believed that this cannot be internalised in ``standard'' homotopy type theory (i.e.\ we cannot write down a type of such diagrams in Agda, Coq or Lean).
 However, we expect that what we do can be emulated in HTS-style systems (such as~\cite{voe_HTS,altCapKra_twolevels}),
 and models in which it is possible have been considered in~\cite{shulman_inversediagrams}.

 For an inverse category, it is easy to see how to get a homotopy coherent diagram with identities from one without identities and vice versa.
 We will sketch a proof that the two types of diagrams are equivalent
 (\cref{subsec:neutralityOfIds}), but
 this again is fairly involved and the complete arguments will be presented in a future full version of the current paper~\cite{we_futurework}.

\paragraph*{\textbf{Diagrams over certain Reedy categories}} 
 
 \emph{Reedy categories} are a generalisation of inverse categories with much weaker conditions on morphisms.
 Prominent examples are the category $\Delta$ (or $\oppo \Delta$), this time not restricted to injective functions, 
 and finite versions of it (i.e.\ $\oppo \Delta$ restricted to objects of level $\leq n$).
 Given a Reedy category, we show how to construct an inverse category with \emph{markings} from it, over which we can then consider Reedy-fibrant diagrams which take the markings into account. 
 The key point is that this inverse category will be finite whenever the Reedy category is,
 allowing us to construct diagrams for a certain class of Reedy categories in a finite manner (which can thus be expressed in ``standard'' homotopy type theory).
 We expect this to be of some interest for the homotopy type theory community as it in particular allows us to add degeneracies to the usual encodings of semisimplicial types.
 In other words, we are able to present a construction of \emph{simplicial types} (with the usual caveat which holds for semisimplicial types, namely that we can only do it ``up to a given level'' in ``standard'' homotopy type theory and require a stronger theory to internalise the full infinite definition).

 For the considered class of Reedy categories (including $\oppo \Delta$), this notion of diagrams can be shown to be equivalent to our notion of general homotopy coherent diagrams.
 Again, the proof is only sketched in \cref{subsec:simplTypesVsHC}.
 In the full version of the current paper~\cite{we_futurework}, we will further present a second construction which is expected to work more generally. 

\subsection*{Organisation}  
We review Reedy fibrant diagrams and explain some constructions with them in \cref{sec:ReedyConstructions}.
In \cref{sec:hCoh}, we propose our notion of homotopy coherent diagrams with a strictification construction, and afterwards (\cref{sec:cohId}), our general homotopy coherent diagrams with identities.
\cref{sec:ReedyDiagsSection} is devoted to the construction of simplicial types and, more generally, diagrams over Reedy categories.
Finally in \cref{sec:equivalencesBetweenConstructions}, we very briefly outline some key ideas belonging to proofs that the different notions of diagrams are equivalent whenever the question makes sense.

\subsection*{Our setting}
Our work should be understood to take place in a \emph{type-theoretic fibration category} with a universe as considered by Shulman~\cite{shulman_inversediagrams}, or a similar setting.
We will present various notions of diagrams, and to ensure that the (large) types of these diagrams are equivalent in the sense of homotopy type theory, we need the universe to be univalent. 
However, the construction of the various diagrams itself, on which we focus here, makes sense without this requirement.

To simplify the presentation, let us pretend that our type-theoretic fibration category is simply the syntactic category of contexts and context morphisms of ``standard'' homotopy type theory (or equivalently the category of types and functions, assuming an $\eta$-law for $\Sigma$-types). 
We write $\type$ for this category. 
In particular, all our examples will be formulated in type theory, and we appeal to type-theoretical intuition.
Thus, familiarity with the terminology of the standard textbook on homotopy type theory~\cite{HoTTbook} is helpful for the examples (and necessary for the proof sketches in \cref{sec:equivalencesBetweenConstructions}), although it is not strictly required to follow the main constructions of the paper (\cref{sec:ReedyConstructions,sec:hCoh,sec:ReedyDiagsSection}).
Familiarity with the basic categorical nerve construction will be helpful.

\section{Constructions with Reedy Fibrant Diagrams} \label{sec:ReedyConstructions}

Although we can in general not encode \emph{strict} type-valued diagrams in type theory, it is possible to do this for certain well-behaved index categories called \emph{inverse categories} and so-called \emph{Reedy fibrant diagrams} over them in an inductive manner.
This is important for us for at least two reasons.
First, we need them as important tools in our constructions.
Second, we can use them as a reference with which we can compare the different kinds of diagrams we will construct. 
The theory of Reedy fibrant diagrams in type theory has been explored in detail by Shulman~\cite{shulman_inversediagrams}, which we very briefly introduce in \cref{sec:introtoreedy}.

\subsection{Inverse categories and Reedy fibrant diagrams} \label{sec:introtoreedy}

Let us begin with the definition of an inverse category.
For the sake of simplicity of the presentation of this short paper, we choose to define a special case.
\begin{definition}[inverse category] \label{def:inverseCat}
 A category $I$ is said to be \emph{inverse} if we can associate with every object a natural number (its \emph{degree}) such that for every non-identity morphism, the degree of its codomain is lower than the degree of its domain.
 We also add the requirement that every object is the domain of a finite number or morphisms.
\end{definition}
 
An alternative description is the following:
Let us write $\N$ for the poset of natural numbers, viewed as a category.
Then, $I$ is an inverse category if we have a functor $F \co I \to \Nop$ that reflects identities, i.e.\ $f$ is an identity whenever $F(f)$ is.
Again, we require every object to be the domain of only finitely many morphisms.

A category $C$ is \emph{direct} if $\oppo C$ is inverse. (Often, \emph{inverse} is defined in terms of \emph{direct}. We have chosen to do the opposite since \emph{inverse} is more central in our development.)

\begin{example} \label{ex:inverse}
As a running example, we will use the following inverse category $\E$ (``example'') with objects $\{x,y,z\}$ and morphisms generated by $u, v \in \E(y, x)$ and $w \in \E(z, y)$ subject to $u \circ w = v \circ w$:
\[
\begin{tikzpicture}
\node (z) {$z$};
\node (y) [right of=z] {$y$} ;
\node (x) [right of=y] {$x$\rlap{.}};

\draw[->] (z) to node [above] {$w$} (y);
\draw[->, bend left =30] (y) to node [above] {$u$} (x);
\draw[->, bend right=30] (y) to node [below] {$v$} (x);
\end{tikzpicture}
\]
We can take the degrees of $x$, $y$, and $z$ to be $0$, $1$, and $2$, respectively. 
\end{example}

The example category $\E$ is non-trivial but sufficiently small to allow the explicit demonstration of some constructions.

The crux of an inverse category $I$ is that certain type-valued diagrams $R$ over $I$ can be defined inductively.
Some preparation is required to make this precise.
Assume $x$ is an object of $I$. Following Shulman~\cite{shulman_inversediagrams}, we write $x \sslash I$ for the full subcategory of the co-slice category $x \slash I$ with the object $\mathsf{id}_x$ removed.
This means that objects of $x \sslash I$ are pairs $(y,f)$ with $y$ of lower degree than $x$ and $f \in I(x,y)$; and morphisms between $(y,f)$ and $(z,g)$ are morphisms $h \in I(y,z)$ such that $h \circ f \equiv g$, as usual for co-slice categories.
$x \sslash I$ is again an inverse category, and we have a canonical forgetful functor $U \co x \sslash I \to I$, mapping $(y,f)$ to $y$.
Because of our additional assumption on inverse categories (see \cref{def:inverseCat}), $x \sslash I$ is always finite, independently of whether $I$ is.
Given a functor $R \co I \to \type$, the \emph{matching object} $M^R_x$ of $R$ at $x$ is defined as the limit of $R \circ U$, 
\begin{equation} \label{eq:matching}
 M^R_x \defeq \lim_{x \sslash I} (R \circ U).
\end{equation}
Of course, limits do a priori not necessarily exist in $\type$, but for the functors we are interested in, this one will exist as shown by Shulman~\cite{shulman_inversediagrams}.
We will see examples in a moment.
The following definition is standard:
\begin{definition}
 A functor $R \co I \to \type$ is called \emph{Reedy fibrant} if, for every object $x$ of $I$, the canonical map $R_x \to M^R_x$ can be written as a projection (or a \emph{fibration}), i.e.\ as a context morphism which simply forgets some of the entries of the context.
\end{definition}
This allows the following inductive construction.
To define a Reedy-fibrant diagram $R$ at an object $x$, assume that we have already defined $R$ on the full subcategory of all objects of lower degree than $x$.
That means in particular that $M^R_x$ is already defined.
To define $R$ at the object $x$, it is then enough to give a type in context $M^R_x$.

\begin{example} \label{ex:strict-diagram}
 In the empty context, a Reedy fibrant diagram over the example category $\E$ from \cref{ex:inverse} is given by three types 
 $R_x$, $R_y$, and $R_z$ in the contexts containing the three matching objects $M^R_x$, $M^R_y$, and $M^R_z$.
 We start with $x$.
 As the category $x \sslash \E$ is empty, we have that $M^R_x$ is the unit type which we can safely ignore and say
 \begin{equation}
  \vdash R_x \istype.
 \end{equation}
 The matching object $M^R_y$ can easily be calculated as $(a : R_x), (b : R_x)$; hence, $R_y$ will be a type in this context,
 \begin{equation}
  (a : R_x), (b : R_x) \vdash R_y \istype.
 \end{equation}
 By calculation, we see that the matching object $M^R_z$ is $(a : R_x), (l : R_y(a,a))$, and we thus need
 \begin{equation}
  (a : R_x), (l : R_y(a,a)) \vdash R_z \istype.
 \end{equation}
 We can think of $R_x$ as a type of \emph{points}.
 For any two points $a,b$, we have a type $R_y(a,b)$ of \emph{lines} between these points.
 Whenever we have a loop, i.e.\ $a:R_x$ and $l : R_y(a,a)$, we have a type $R_z(a,l)$ of \emph{fillers} for this loop.
 From this point of view, we see the matching objects as types of ``boundaries''.

 These three pieces of data encode the strict diagram $R \co \E \to \type$ as follows:
 \begin{itemize}
  \item $x$ is mapped to $(a : R_x)$, a context of length one;
  \item $y$ is mapped to $(a : R_x), (b : R_x), (l : R_y(a,b))$, a context of length three;
  \item and, finally, $z$ is mapped to $(a : R_x), (l : R_y(a,a)), (f : R_z(a,l))$, another context of length three.
 \end{itemize}
 The morphisms of $\E$ are mapped to projections (and duplications) between contexts, e.g.\ $w$ is mapped to $(a,l,f) \mapsto (a,a,l)$.
 As composition of projections is strictly associative, the functor laws automatically hold strictly.
 Therefore, there is no need to include them explicitly in the definition of $R$, which is the whole trick of this presentation.
 
 More sophisticated inverse categories could encode more complicated structure and for example contain a type of fillers for triangles instead of only for loops.
\end{example}

When $I$ is finite, we can easily internalise the construction, i.e.\ perform it completely inside type theory.
One can then talk about a \emph{type of $I$-diagrams}.
This is something one may want to do in order to implement it in a proof assistant.

\begin{example} \label{ex:E-Reedyfibrant-internally}
 Let $\UU$ be a type universe.
 Inside type theory, the type of Reedy fibrant diagrams over $\E$ valued in $\UU$ is a (nested) $\Sigma$-type with the following three components:
 \begin{align*}
  &R_x : \UU \\
  &R_y : R_x \times R_x \to \UU \\
  &R_z : \left(\Sigma (a:R_x). R_y(x,x)\right) \to \UU.
 \end{align*}
\end{example}

\subsection{The simplex category and semisimplicial types} \label{sec:DeltaAndSemi}

Recall the definition of the category $\Delta$, which is particularly important in algebraic topology.
It has non-empty finite ordinals as objects, written as $[0], [1], [2], \ldots$ where $[n]$ is the ordinal $\{0,1,\ldots,n\}$.
The morphisms from $[m]$ to $[n]$ consist of monotone (order-preserving, but not necessarily injective) functions $[m] \to [n]$.
There is a canonical degree map, $[m] \mapsto m$, but note that morphisms can both increase or decrease this degree.
To remedy this, we can consider the subcategory $\deltplus$ of $\Delta$ which has only strictly monotone (injective order-preserving) functions as morphisms.
This is a direct category; its opposite $\deltop$ is an example of an inverse category and very important to us.
A Reedy fibrant diagram over $\deltop$ is known in the homotopy type theory community as a \emph{semisimplicial type}~\cite{herbelin_semisimpl,nicolai:thesis,shulman:eating}.
Up to level $2$, such a semisimplicial type consists of a type of points, a type of lines, and a type of triangles, as in:
\begin{align*}
 &\vdash A_{[0]} \istype \\[0.3em]
 (p_0 : A_{[0]}), (p_1 : A_{[0]}) &\vdash A_{[1]} \istype \\[0.3em]
 (p_0 : A_{[0]}), (p_1 : A_{[0]}), (p_2 : A_{[0]}), (l_{01} : A_{[1]}(p_0,p_1)), \span \\
   (l_{12} : A_{[1]}(p_1,p_2)), (l_{02} : A_{[2]}(p_0,p_1)) &\vdash A_{[2]} \istype
\end{align*}
As in \cref{ex:E-Reedyfibrant-internally}, we can easily present this as a collection of type families in type theory as long as we only want to do it up to a fixed finite level $n$ (see e.g.~\cite{kraus_haskell} for a Haskell script generating the relevant Agda code). 
Whether the full definition with components $A_{[n]}$ for all $n$ can be encoded in ``standard'' homotopy type theory is the well-known open problem mentioned in the introduction, and the development of HTS-style systems~\cite{voe_HTS,altCapKra_twolevels,DBLP:journals/corr/PartL15} has been inspired by the desire to perform this construction.

\subsection{A Reedy fibrant replacement construction}

Given an inverse category $I$, we can construct a Reedy fibrant functor $R \co I \to \type$ by induction on the objects of $I$.
Another possibility is to start with a (non-fibrant) functor that is given in another way, and somehow build a Reedy fibrant version of it, a \emph{Reedy fibrant replacement}.

One possibility is as follows.
Assume that $F \co I \to \type$ is a functor, i.e.\ every $F_i$ is in a context, every $F(f)$ is a function between contexts, and composition of $F(g)$ and $F(f)$ happens to be strictly (i.e.\ judgmentally) equal to $F(g \circ f)$ whenever $g$ and $f$ compose in $I$.

In particular, we are interested in the strict functor $\Sp \co \deltop \to \type$ which is given by
\begin{align*}
 \Sp_{[n]} \; \defeq \; &(A_0 : \UU), (A_1 : \UU), (A_2 : \UU), \ldots, (A_n : \UU), \\
                        &(f_0 : A_0 \to A_1), (f_1 : A_1 \to A_2), \ldots, \\
                        &(f_{n-1} : A_{n-1} \to A_n).
\end{align*}
This is the nerve of the internal category $\UU$, restricted from $\Delta$ to $\deltplus$.
In short, $\Sp_{[n]}$ is a context of $(n+1)$ types and $n$ functions between them.
The morphism part of $\Sp$ is given in the canonical way: for an injective monotone function $g \co [m] \to [n]$, the function $\Sp(g) : \Sp_{[n]} \to \Sp_{[m]}$ throws away some of the types and composes the functions accordingly.
Using that function composition is strictly associative, this gives rise to a functor.

We can now construct a Reedy fibrant functor $G \co I \to \type$ together with a natural transformation $\eta \co F \to G$ which is levelwise a homotopy equivalence (which means that, if we regard contexts as nested $\Sigma$-types, then each map $\eta_x$ is an equivalence between the type $F_x$ and $G_x$). 
This is only a special case of the more general construction by Shulman~\cite{shulman_inversediagrams} who uses that every natural transformation, here the unique one from $G$ to the terminal functor, factors as an \emph{acyclic cofibration} followed by a Reedy fibration.
It has been given explicitly for the case that $G$ is a constant functor in~\cite{kraus_generaluniversalproperty}, and the case we are interested in has been suggested in~\cite{altenCapKraus_infCats}.

A concrete construction of $G$ and $\eta \co F \to G$ can be done as follows, by induction on $I$.
Say, $x$ is an object of $I$ and both $G$ and $\eta$ are defined for all objects of smaller degree than $x$.
Then, $M^G_x$ is fully specified, and $\eta$ extends to a function $\tilde \eta \co F_x \to M^G_x$.
In context $m : M^G_x$, we define $G_x$ to be $\Sigma (a : F_x). (\tilde \eta(a) = m)$, 
and we can extend $\eta$ to $x$ using the map $F_x \to G_x$ sending $a$ to $(\tilde \eta(a), a, \refl)$.
Note that the pair $(\tilde \eta(a), \refl)$ is a ``singleton'' which is known to inhabit a contractible type, ensuring that $\eta_x$ is an equivalence.

The point of the strict functor $\Sp$ is that its fibrant replacement makes precise the idea of having a semisimplicial type $\ty$ such that 
$\ty_{[0]}$ are types (i.e.\ $\ty_{[0]}$ is $\UU$), for two types $A,B$, we have that $\ty_{[1]}(A,B)$ is the type of functions $A \to B$, 
for three types and functions $A$, $B$, $C$, $f : A \to B$, $g : B \to C$, $h : A \to C$, 
we have that $\ty_{[2]}(A,B,C,f,g,h)$ is the equality type $g \circ f = h$, and so on.
Strictly speaking, if we take $\ty$ to be the fibrant replacement of $\Sp$ as outlined above, $T_{[0]}$ will not \emph{judgmentally} be $\UU$, 
but only equivalent to it, and similarly for $\ty_{[1]}$ and $\ty_{[2]}$.
However, we can always manually tweak $\ty$ on a finite number of levels, 
and for the sake of a nicer presentation, we therefore assume that $T$ is really types, 
functions, and proofs of commutativity on the lowest three levels. 

\begin{remark} \label{restricted-shulman-universe}
We can get a different model for the fibrant replacement of $\Sp$ by modifying Shulman's universe~\cite{shulman_inversediagrams} for inverse diagrams over $\deltplus$.
Shulman's universe $V$ can be seen as a semisimplicial type of types, relations between types, higher relations between a triangle of relations, etc.
By adding appropriate propositional constraints to its components, we can change it into a semisimplicial type of types, functions between types, and the higher behaviour we are looking for.
Note that this makes sense even without the $\eta$-law for functions, which the definition of $\Sp$ depended upon.
\end{remark}

\subsection{Downwards closed full subcategories and Reedy limits} \label{sec:DownwardsClosedIntuition}

If we have a finite inverse category $I$ and a Reedy fibrant functor $R \co I \to \type$, we can form its limit $\lim_{I} R$, which will be a context in $R$.
We have already seen this in the construction of the matching object, and it is true in general that this limit exists in $\type$~\cite{shulman_inversediagrams}.
Intuitively, we think of this limit as a context containing one component for every object in $I$.
We can easily turn such a finite context into a type: we just form a nested $\Sigma$-type (or a \emph{record}, if we allow our type theory to have that notion) with one component for every entry of the context.

Now assume that $J$ is a ``downwards closed full subcategory'' of $I$.
What we mean by this is that $J$ is a full subcategory $J \subset I$ with the property that, if $x$ is an object in $J$ and there is a morphism $f \in I(x,y)$, then $y$ (and automatically $f$) is also present in $J$.
We can form the limit of $R$ restricted to $J$, written $\lim_{J} R$.
A very intuitive but no less important observation is that this limit will simply be a subcontext of the context $\lim_I R$, 
and the canonical context morphism $\lim_I R \fib \lim_J R$ is a projection (in technical terms, a \emph{fibration} in $\type$) which simply removes some entries of the bigger context.
This can also be seen as a (generalised) projection map between the corresponding nested $\Sigma$-types.

\section{Homotopy coherent diagrams} \label{sec:hCoh}

In this section, we make precise the construction of diagrams which include all higher coherences explicitly.
For the whole section, let us assume that $I$ is some inverse category.
We mostly think of the case that $I$ is finite apart from \cref{sec:cohId} where this condition is dropped.

\subsection{Preliminary Observations on the Positive Nerve} 

Given $I$, the well-known nerve construction yields a simplicial set $\fullnerv I$.
Recall that $n$-cells (elements of $(\fullnerv I)_n$) are given as $n$-strings of composable morphisms $X_0 \xrightarrow{f_1} X_1 \xrightarrow{f_2} \ldots \xrightarrow{f_n} X_n$, for simplicity written $\xrightarrow{f_1}\xrightarrow{f_2}\ldots\xrightarrow{f_n}$.
Not needing the degeneracy structure, we view $\fullnerv I$ as a \emph{semi}simplicial set.
The fact that $I$ is inverse implies that the composition of non-identity arrows is a non-identity arrow.
Thus, we can consider the semisimplicial set whose $n$-cells are $n$-strings of composable non-identity arrows. 
We write $\posnerv I$ for this semisimplicial set and call it the \emph{positive nerve} of $I$.

Further, we can form the category of elements of $\posnerv I$, written $\int \posnerv I$.
If we spell it out, we see that objects of the category $\int \posnerv I$ are sequences of composable non-identity arrows of $I$.
We have a morphism from $\xrightarrow{f_1} \xrightarrow{f_2} \ldots \xrightarrow{f_n}$ to $\xrightarrow{g_1} \xrightarrow{g_2} \ldots \xrightarrow{g_k}$
if the latter sequence can be constructed from the first by composing arrows and by discarding arrows in the beginning and end.

\begin{example} \label{ex:posnervE}
Let us discuss the example $\E$. 
The category $\int \posnerv \E$ has nine objects and can be pictured as shown below, where we denote sequences of length $0$ simply by their single object, and longer ones with their morphisms.
Note that $u \circ w = v \circ w$.
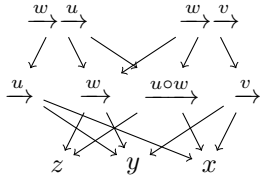
\[
\begin{tikzpicture}[scale=0.5]
\node (wu) at (1,4) {$\xrightarrow{w}\xrightarrow{u}$}; 
\node (wv) at (5,4) {$\xrightarrow{w}\xrightarrow{v}$}; 
\node (u) at (0,2) {$\xrightarrow{u}$}; 
\node (w) at (2,2) {$\xrightarrow{w}$}; 
\node (uw) at (4,2) {$\xrightarrow{u \circ w}$}; 
\node (v) at (6,2) {$\xrightarrow{v}$}; 
\node (z) at (1,0) {$z$};
\node (y) at (3,0) {$y$};
\node (x) at (5,0) {$x$};

\draw[->] (wu) to node {} (w);
\draw[->] (wu) to node {} (u);
\draw[->] (wu) to node {} (uw);
\draw[->] (wv) to node {} (uw);
\draw[->] (wv) to node {} (w);
\draw[->] (wv) to node {} (v);
\draw[->] (u) to node {} (x);
\draw[->] (u) to node {} (y);
\draw[->] (v) to node {} (x);
\draw[->] (v) to node {} (y);
\draw[->] (w) to node {} (y);
\draw[->] (w) to node {} (z);
\draw[->] (uw) to node {} (x);
\draw[->] (uw) to node {} (z);
\end{tikzpicture}
\]
\end{example}

There are several observations to make, heavily using the fact that $I$ is inverse. 
First, we see that $\int \posnerv I$ is inverse again, with the same height as $I$; the degree of a sequence is given by its length.
Second, $\int \posnerv I$ is a preorder (there is at most one morphism between two given objects).
Third, if $I$ is finite, then so is $\int \posnerv I$. 
Fourth, we have a canonical functor $\fst \co \int \posnerv I \to \deltop$, sending a sequence $\xrightarrow{f_1} \ldots \xrightarrow{f_n}$ to $[n]$.
Fifth, assume that $R \co I \to \type$ is Reedy fibrant. Given any functor $F \co A \to I$, it is not in general the case that $R \circ F$ is Reedy fibrant; indeed, many counterexamples can readily be given by taking $I$ to be the terminal category.
However, it is the case that $R \circ \fst$ is always Reedy fibrant (the abstract reason is that $\fst$ is a discrete Grothendieck opfibration).

\subsection{The definition of homotopy coherent diagrams}

For homotopy coherent diagrams, it turns out to be easier to already define it as a notion internal to type theory:
\begin{definition} \label{def:hcoh}
 We define the type of \emph{homotopy coherent diagrams} over $I$ to be the nested $\Sigma$-type corresponding to the limit of the composition
 \begin{equation*}
  \big({\textstyle{\int}} \posnerv I\big) \xrightarrow {\fst} \deltop \xrightarrow{\ty} \type.
 \end{equation*}
\end{definition}
A word of explanation: the limit $\lim_{\int\posnerv I}(\ty \circ \fst)$ exists thanks to the assumed finiteness of $I$ and thanks to $\ty \circ \fst$ being Reedy fibrant. 
This limit is an object in $\type$, hence a context, which can be turned into a nested $\Sigma$-type, and that type is what we call \emph{homotopy coherent} or \emph{weak diagrams}.
If we really want an external notion, we can of course say that a homotopy coherent diagram in context $\Gamma$ is a morphism in $\type$ from $\Gamma$ to the stated limit.
That is, intuitively, a homotopy coherent diagram over $I$ consists of 
\begin{itemize}
 \item a type $X_i$ for every object $i$ of $I$;
 \item a function $X_f : X_i \to X_j$ for every morphism $f \in I(i,j)$;
 \item an equality $X_g \circ X_f = X_{g \circ f}$ for every pair of composable morphisms in $I$;
 \item proofs that these equalities are associative (one associativity proof for every triple of composable morphisms);
 \item all the higher dimensional associahedra and so on, expressing higher coherence laws.
\end{itemize}
To continue with our running example, let us see how this works out for $\E$:
\begin{example} \label{ex:weak-diagram}
A homotopy coherent diagram over $\E$ has nine components, corresponding to the nine objects of $\int \posnerv \E$ (see \cref{ex:posnervE}):
\begin{itemize}
 \item three types $X$, $Y$, and $Z$;
 \item three functions: $\ww : Z \to Y$ and $\uu,\vv : Y \to X$, reusing the names of the morphisms in $\E$;
 \item a further function $\sss : Z \to X$, corresponding to the morphism $u \circ w$ (which is also $v \circ w$);
 \item and two equalities, $\pp : \sss = \uu \circ \ww$ and $\qq : \sss = \vv \circ \ww$.
\end{itemize}
\end{example}

\subsection{From a strict diagram to a weak diagram}
Given a Reedy fibrant diagram $A$ over $I$ (with values in $\UU$), we want to construct a homotopy coherent diagram.
Writing $\unit$ for the terminal object of $\type$, we need to define a cone $\mathsf{const}_\unit \to (\ty \circ \fst)$.
Objects of $\int \posnerv I$ are of the form $i_0 \xrightarrow{f_0} i_1 \xrightarrow{f_1} \ldots i_n$.
A cone $\mathsf{const}_\unit \to (\Sp \circ \fst)$ is given by choosing the component at an object of this form to be $A(i_0) \xrightarrow{A(f_0)} A(i_1) \xrightarrow{A(f_1)} \ldots A(i_n)$, formally the context $(A(i_0) : \UU), \ldots, (A(i_n) : \UU), (A(f_0) : A(i_0) \to A(i_1)), \ldots$.
We then compose with $\eta \co \Sp \to \ty$.

To internalise this construction (as always for a fixed index category $I$),
we replace $\unit$ by the type encoding Reedy fibrant diagrams over $I$, 
and we replace all occurrences of $A$ by the corresponding projections.
This gives a function in type theory which turns an element of the type of strict diagrams into an element of the type of weak diagrams.

\begin{example}
 Assume we are given a Reedy fibrant diagram over $\E$ as in \cref{ex:E-Reedyfibrant-internally}, i.e.\ we are given $R_x$, $R_y$, $R_z$.
 We want to calculate the corresponding weak diagram as in \cref{ex:weak-diagram}.
 Recall that such a weak diagram first of all consists of three types $X$ and $Y$ and $Z$, which are here given as $R_x$ and $\Sigma (a,b : R_x). R_y(a,b)$ and $\Sigma (a:R_x). \Sigma (l : R_y(a,a)). R_z(a,l)$.
 The function $\ww : Z \to Y$ is given by $(a,l,f) \mapsto (a,a,l)$, the function $\vv : Y \to X$ is given by $(a,b,l) \mapsto b$, and so on.
 The required equalities hold on the nose, $\pp$ and $\qq$ are just $\refl$.
\end{example}

\subsection{A strictification construction} \label{sec:strictification}
Naturally, the construction of a strict diagram from a weak one is more involved.
We write $I + \iota$ for the category $I$ with one additional object $\iota$ formally added. 
Note that this is really the coproduct of $I$ and the terminal category, $I + \mathsf{1}$; we call the added object $\iota$ simply to have a name to refer to it. 
For a given object $i$ of $I$, we further write $I + \iota_{\to i}$ for the category $I + \iota$ with a morphism $\iota \to i$ added (which of course freely generates a morphism $\iota \to j$ for every morphism $i \to j$ in $I$).
Then, $\int \posnerv (I + \iota)$ is a ``downwards closed full subcategory'' 
of $\int \posnerv(I + \iota_{\to i})$ (see \cref{sec:DownwardsClosedIntuition}), with the additional objects being
all sequences of non-identities of the form $\iota \to i_0 \xrightarrow{f_0} i_1 \xrightarrow {f_1} \ldots \xrightarrow{f_{n-1}} i_n$, with $n \geq 0$, and where $i_0$ could be $i$.

Let $X$ be a homotopy coherent diagram over $I$.
This gives rise to an arrow $\unit \xrightarrow{X,\unit} \lim_{\int \posnerv (I + \iota)}(\ty \circ \fst)$ by constructing the corresponding cone, where the single new component of the zero-length sequence $\iota$ is given by the unit type $\unit$, and the rest by $X$.

To define the strict diagram $A \co I \to \type$ we form, for any object $i$, the following pullback:
\[
 \begin{tikzpicture}[x=4cm,y=-1.5cm]
  \node (A) at (0,0) {$A_i$}; 
  \node (U) at (0,1) {$\unit$}; 
  \node (Iiarrow) at (1,0) {$\lim_{\int\posnerv(I + \iota_{\to i})}(\ty \circ \fst)$}; 
  \node (Ii) at (1,1) {$\lim_{\int\posnerv(I + \iota)}(\ty \circ \fst)$};

  \draw[->>,dashed] (A) to node [] {} (U);
  \draw[->>] (Iiarrow) to node [] {} (Ii);
  \draw[->,dashed] (A) to node [] {} (Iiarrow);
  \draw[->] (U) to node [above] {$(X, \unit)$} (Ii);
 \end{tikzpicture}
\]
In other (more type-theoretic) words, the type $A_i$ is given as the fibre over the element of $\lim_{\int\posnerv(I + \iota)}(\ty \circ \fst)$ that is given by $(X,\unit)$; it can be thought of as a nested $\Sigma$-type with one component for each sequence of positive length starting with $\iota \to \ldots$ (see \cref{ex:strictification} below).
For a morphism $f \in I(i,j)$, we have $A(f) : A_i \to A_j$ given by projection, as $\int \posnerv(I + \iota_{\to j})$ will be a ``downwards closed full subcategory'' of $\int \posnerv(I + \iota_{\to i})$, using that $I$ is inverse.

Finally, we need to check that the such-defined strict diagram $A \co I \to \type$ is Reedy fibrant.
For an object $i$, let us write $I + \iota_{\dashrightarrow i}$ for the category $I + \iota_{\to i}$, with the single arrow $\iota \to i$ removed (but keeping all the morphisms generated by it); the fact that removing this single arrow makes sense uses once more that $I$ is inverse.
We claim that the matching object $M^A_i$ is given as the pullback of $(X,\unit)$ along the projection corresponding to the inclusion of categories $\int \posnerv(I + \iota) \; \subset \; \int \posnerv(I + \iota_{\dashrightarrow i})$.
After we verify this, we are done, as the inclusion $\int \posnerv \left(I + \iota_{\dashrightarrow i}\right) \; \subset \; \int \posnerv \left(I + \iota_{\to i}\right)$ gives rise to a projection $A_i \fib M^A_i$. 
In other words, the components of $A_i$ not present in $M^A_i$ are those corresponding to sequences starting with $\iota \to i \to \ldots$.

To verify the claim, we perform a standard calculation as follows.
By definition and rewriting, we have
\begin{equation}
 M^A_i \; \cong \; \lim_{i \sslash I} (A \circ \fst) \; \cong \; \lim_{(x,f) \in i \sslash I} A_x
\end{equation}
As $A_x$ is itself defined as a limit (a pullback), we can commute these limits.
Two of the objects in the cospan defining $A_x$ are independent of $x$ (there is nothing to do when taking their limit), 
and we get that the object above is isomorphic to the pullback of
\[
 \begin{tikzpicture}[x=4cm,y=-1.5cm]
  \node (U) at (0,1) {$\unit$}; 
  \node (Iiarrow) at (1,0) {$\lim_{(x,f) \in i \sslash I} \left(\lim_{\int\posnerv(I + \iota_{\to x})}(\ty \circ \fst) \right)$}; 
  \node (Ii) at (1,1) {$\lim_{\int\posnerv(I + \iota)}(\ty \circ \fst)$}; 

  \draw[->>] (Iiarrow) to node [] {} (Ii);
  \draw[->] (U) to node [above] {$(X, \unit)$} (Ii);
 \end{tikzpicture}
\]

Combining the index categories of the ``nested limit'' gives $\int \posnerv(I + \iota_{\dashrightarrow i})$ as claimed.

\begin{example} \label{ex:strictification}
 To continue with our running example, let us assume we are given a weak diagram over $\E$.
 As in \cref{ex:weak-diagram}, we assume we are given this diagram as three types $X$, $Y$, $Z$, corresponding to the objects of $\E$; functions $\ww$, $\uu$, $\vv$, corresponding to the morphisms of $\E$; one more function $\sss : Z \to X$; and two equalities $\pp$, $\qq$, expressing the connection between $\uu$, $\vv$, $\ww$, and $\sss$.
 We wish to construct a Reedy fibrant diagram $A$ from this. 
 Recall that $\E$ has objects $x$, $y$, $z$.
 Compared to $\int \posnerv \left(I + \iota\right)$, we see that:
 \begin{itemize}
  \item $\int \posnerv \left(I + \iota_{\to x}\right)$ has exactly one additional object, namely $\iota \to x$.
  \item $\int \posnerv \left(I + \iota_{\to y}\right)$ has an additional objects $\iota \to y$ (let us call this $\eta$), but also $\iota \xrightarrow{u \circ \eta} x$ and $\iota \xrightarrow{v \circ \eta} x$ and, finally, $\iota \to y \xrightarrow{u} x$ and $\iota \to y \xrightarrow{v} x$.
  \item $\int \posnerv \left(I + \iota_{\to z}\right)$ has one object $\iota \to x$, one object $\iota \to y$, one object $\iota \to z$, two objects $\iota \to y \to x$, one object $\iota \to z \to y$, one object $\iota \to z \to x$, and two objects $\iota \to z \to y \to x$.
 \end{itemize}
 To get an explicit listing of the components of the matching object $M^A_\epsilon$, we simply need to take $A_\epsilon$ and remove all sequences containing $\epsilon$ (where $\epsilon \in \{x,y,z\}$). 
 
 The Reedy fibrant diagram that we get by using the above formula can be represented as follows, giving $A_\epsilon$ as contexts over the relevant matching object $M^A_\epsilon$: 
 \begin{align*}
  A_x \defeq  \, &X \\
  M^A_y \defeq \, & (x_1 : X), (x_2 : X) \\
  A_y(x_1,x_2) \defeq \, &(y : Y), \\
    &(\delta_\uu : \uu(y) = x_1), (\delta_\vv : \vv(y) = x_2) \\
  M^A_z \defeq \, &(x : X), (y : Y), \\
    &(\delta_\uu : \uu(y) = x), (\delta_\vv : \vv(y) = x) \\
  A_z(x,y,\delta_\uu,\delta_\vv) \defeq \, &(z : Z), \\
    &(\omega: \ww(z) = y), (\sigma : \sss(z) = x), \\
    &(\Theta_\uu : \sigma = \pp(z) \ct \mapfunc_\uu(\omega) \ct \delta_\uu), \\
    &(\Theta_\vv : \sigma = \qq(z) \ct \mapfunc_\vv(\omega) \ct \delta_\vv)
 \end{align*}
 Note that we have $\pp : \sss = \uu \circ \ww$ and we write $\pp(z)$ for the equality we get by applying both sides to $z$.
 This is sometimes written as $\mathsf{happly}_\pp(z)$ in homotopy type theory.
\end{example}
The strictification construction can be internalised to a function from the type of weak diagrams to the type of Reedy fibrant diagrams over $I$.
We will show in~\cite{we_futurework} that this function is an equivalence, a sketch of which can be found in \cref{subsec:weakStrictEquiv}.

\subsection{Homotopy coherent diagrams with identities} \label{sec:cohId} 

Our goal of this section is to define the notion of a \emph{general homotopy coherent diagram} which works for index categories which might not be inverse.

So far, we have heavily used a simple but powerful fact: as long as we restrict ourselves to \emph{inverse} categories $I$, we can completely ignore identity morphisms (recall that inverse categories are in one-to-one correspondence to inverse semicategories).
Together with the assumption that $I$ is finite, this implies the nice property that $\int \posnerv I$ is finite.
We know that any Reedy fibrant diagram over a finite index category has a limit~\cite{shulman_inversediagrams}, which has allowed us to perform all described constructions in ``standard'' homotopy type theory.

The situation is more involved if we want to go further and consider an index category $C$ which is not necessarily inverse.
In this case, we have to take identities into account as the composition of two non-identities may very well be an identity (in other words, many categories are \emph{not} freely generated from semicategories).
We now will do this and consider the category of elements of the full nerve, $\int \fullnerv C$.
Note that $\int \fullnerv C$ is still an inverse category, and in fact quite well-behaved (it fulfills the condition mentioned in \cref{def:inverseCat}, i.e.\ every coslice of it is finite).
However, as long as $C$ has at least one object, it will always be infinite.
This means we can weaken our previous assumption that the index categories are finite, as it now does not make a difference anymore; a canonical example which is of interest is the category $\Delta^{\mathsf{op}}$. 

While we can externally consider Reedy fibrant diagrams over the infinite index category $\int \fullnerv C$ without problems,
we cannot internalise the construction in ``standard'' homotopy type theory as we would need some sort of ``infinitely nested $\Sigma$-types'', or, in the more precise terms of~\cite{shulman_inversediagrams}, \emph{Reedy $\omega^{\mathsf{op}}$-limits}.
Such \emph{infinitary type theories} have been considered before~\cite{shulman_inversediagrams,kraus_generaluniversalproperty} and are supported by many models.
Further, we expect that we are able to do sufficient encodings in two-level systems such as~\cite{voe_HTS,altCapKra_twolevels}.
For internalisability, let us thus assume for the current section that we work with a type theory where such infinitary constructions are available.

Even with this assumption, it is not so easy to say what a \emph{general homotopy coherent diagram} over $C$ is (we could also call it a \emph{homotopy coherent diagram with identities}).
Let us for a moment go back to an inverse category $I$.
The category $\int \fullnerv I$ has $\int \posnerv I$ as a full subcategory, and the additional objects are those sequences $i_0 \xrightarrow{f_0} \ldots \xrightarrow {f_{n-1}} i_n$ which contain at least one identity morphism.
Thus, if we consider an element $h$ of the limit $\lim_{\int \fullnerv I} (\ty \circ \fst)$, then $h$ will be an ``infinite tuple'' which first of all contains the same components as an ordinary weak diagram over $I$, but in addition $h$ will contain one component $h(s)$ for each sequence $s$ containing at least one identity.
This is not yet what we want.
For an object $x$, the identity $\id_x \in I(x,x)$ can be seen as a sequence of length one, and we need to ensure that $h(\id_x) : h(x) \to h(x)$ is not \emph{any} function, but the actual identity function.
Similarly, longer sequences containing identities need to give rise to actual degeneracies (trivial proofs, i.e.\ ``reflexivities''), not just \emph{any} proofs.
A plausible approach would be to add an equality stating that $h(\id_x)$ is equal to the identity function, and so on.
However, this would then require coherence laws ensuring that these equalities fit together on higher levels, and it is unclear to us whether this approach would ultimately be feasible.

Fortunately, there is a more elegant solution. 
A related setting (outside of type theory) is the following.
A \emph{semi-Segal space} is a diagram over $\deltop$ valued in spaces which corresponds to our semisimplicial types with some additional conditions.
Given such a semi-Segal space, it is natural to ask whether a Segal space can be constructed from it (i.e.\ whether degeneracies can be added).
A very minimalistic strategy was suggested by Harpaz~\cite{harpaz2015quasi}, based on the observation that, on the lowest level, it is sufficient to require enough \emph{equivalences} to exist rather than identities, and that those equivalences can be used to generate the whole degeneracy structure.

It was pointed out to us by Paolo Capriotti that this trick can also be applied in type theory.
The advantage is that, if we add the property that a certain function is an equivalence, we only add a mere proposition which will not require further coherence laws.
It turns out that the following very minimalistic modification of the definition of weak diagrams is sufficient to derive a notion of general homotopy coherent diagrams with identities:

\begin{definition} \label{def:hcohWithId}
 For a category $C$, a \emph{general homotopy coherent diagram} is an element
 \begin{equation} \label{eq:typeOfHC}
  h : \lim_{\int \fullnerv C} (\ty \circ \fst),
 \end{equation}
 of which we think as a large nested tuple with one component for each sequence of arrows in $C$, with the following condition:
 for each object $x$ of $C$, 
 the function $h(\id_x) : h(x) \to h(x)$ is an equivalence.
\end{definition}

\begin{remark}[Relation to Szumilo's D construction]
If we construct the functor $\ty$ as outlined in \cref{restricted-shulman-universe}, we can give a different description of (general) homotopy coherent diagrams.
A \emph{relative category} is a category with a wide subcategory of morphisms called \emph{marked}.
We can generalise our type of Reedy fibrant diagrams over an inverse category to a type of Reedy fibrant diagrams over a relative inverse category where the marked morphisms get mapped to equivalences.

Given an inverse category $I$, we have a functor $\operatorname{fst} \co \int \posnerv I \to I$ given by the first vertex projection.
We let $\operatorname{Sd}_I$ be the relative category given by $\int \posnerv I$ with markings created by $\operatorname{fst}$ (viewing $I$ as a relative category with only identities marked).
Then homotopy coherent diagrams over $I$ are just Reedy fibrant diagrams over $\operatorname{Sd}_I$.

Similarly, given a category $C$, we have a functor $\operatorname{fst} \co \int \fullnerv C \to C$ given by the first vertex projection.
Szumilo's \emph{D construction} $D_C$~\cite{szumilo:thesis} is the relative category given by $\int \fullnerv C$ with markings created by $\fst$.
Now general homotopy coherent diagrams over $C$ correspond to Reedy fibrant diagrams over $D_C$. 

Note that the preceding two paragraphs immediately generalize to relative (inverse) categories $I$ and $C$. 
\end{remark}

We think that the notion of general homotopy coherent diagrams as given in \cref{def:hcohWithId} is less intuitive than the notions of Reedy fibrant and weak diagrams over inverse categories, in the same way as we find Harpaz' result surprising.
Some evidence for our claim that \cref{def:hcohWithId} is well-behaved is given by the following.

First, note that, for an inverse category $I$, we can very easily construct a weak diagram (without identities) from a homotopy coherent one with identities: it is essentially given by a projection which simply removes all components that belong to sequences containing identities.
In \cref{subsec:neutralityOfIds}, we will sketch a proof that this projection is an equivalence.
In other words, for inverse categories, the notion of homotopy coherent diagrams with identities coincides with the ones without identities, and hence also with the Reedy fibrant ones.
(Note that it is not hard to construct an inverse to the mentioned projection explicitly, although we do not need this. 
If we are given a weak diagram and want to add compoents for sequences with identities, we simply add the actual identity functions and, by induction, one sees that one can add trivial components on all higher levels.)

\section{Diagrams Over certain Reedy Categories} \label{sec:ReedyDiagsSection}

A Reedy category $R$ is a category where every object has a degree, just as an inverse category.
The difference is that there is no restriction on the direction in which morphisms can go.
However, a condition is that there are two subcategories of $R$, denoted by $\plus R$ and $\minus R$, such that:
\begin{itemize}
 \item both $\plus R$ and $\minus R$ contain all objects of $R$ (``wide subcategories''),
 \item $\minus R$ is inverse and $\plus R$ is direct (recall that this just means that $\opandplus R$ is inverse), such that $\minus R$ and $\plus R$ use the same degree function,
 \item every morphism of $R$ factors uniquely as a morphism in $\minus R$ followed by a morphism in $\plus R$. 
\end{itemize}

A standard example is the category $\Delta$.
We have introduced $\Delta$ in \cref{sec:DeltaAndSemi}, but immediately restricted to $\deltplus$, containing only \emph{injective} maps.
If we write $\minus \Delta$ for the collection of all \emph{surjective} maps of $\Delta$, we can check that $\Delta$ is indeed a Reedy category, and so is its opposite $\oppo \Delta$ (with $+$ and $-$ switched).

In this section, we want to present a construction of diagrams over certain Reedy categories. 
Compared to our homotopy coherent diagrams with identities, the advantage is that we do no longer require ``infinitary'' constructions; that is, for a finite Reedy category, the construction fully works in ``standard'' homotopy type theory.

The strategy is to replace $R$ by a suitable direct category $D(R)$ with \emph{markings}.
The desired diagrams are then Reedy fibrant diagrams over $\oppo{D(R)}$ with the additional condition that every projection corresponding to a marked morphism is an equivalence.
Again, as we add a merely propositional property, we can avoid the need for further coherence laws.

Before presenting the general construction, we look at the special case of $\Delta$, which allows us to construct \emph{simplicial types}.

\subsection{Simplicial Types}

Recall that we identify an object $[n]$ of $\Delta$ with the preorder $\braces{0, 1, \ldots, n}$.
We first define a \emph{direct replacement} $\D$ of $\Delta$.
We will then be interested in certain Reedy-fibrant diagrams over $\oppo{\D}$.

\begin{definition}
 The category $\D$ is defined as follows.
 Objects are non-empty lists of positive integers, written as $(a_0, a_1, \ldots, a_m)$.
 Morphisms from $(a_0, \ldots, a_m)$ to $(b_0, \ldots, b_n)$ are those morphisms $f \in \Delta([m],[n])$ such that for all $j \in [n]$, we have that $b_j$ is at least as large as the sum of all $a_i$ with $f(i) = j$:
 \begin{align*}
  \D(&(a_0,\ldots,a_m) , (b_0, \ldots, b_n)) \; \defeq \; \\
     &\set[\Big]{f \in \Delta([m],[n])}{\forall j \in [n], b_j \geq \sum_{f(i) = j} a_i}
 \end{align*}
 Composition of morphisms in $\D$ is defined in the canonical way as in $\Delta$.
 We say that a morphism in $\D$ is \emph{marked} if it comes from an identity in $\Delta$.
\end{definition}

\begin{figure}[h]  
\[
\begin{tikzpicture}[x=2cm,y=-1.5cm]
\node(a) at (0,0) {$(1)$};
\node(aa) at (1,0) {$(1,1)$};
\node(aaa) at (2,0) {$(1,1,1)$};
\node(b) at (0,1) {$(2)$};
\node(ba) at (1,1) {$(2,1)$};
\node(ab) at (1.3,1.3) {$(1,2)$};
\node(c) at (0,2) {$(3)$};

\draw[->, transform canvas={yshift=-0.25ex}] (a) to node {} (aa);
\draw[->, transform canvas={yshift=+0.25ex}] (a) to node {} (aa);

\draw[->, transform canvas={yshift=-0.5ex}] (aa) to node {} (aaa);
\draw[->, transform canvas={yshift=+0.0ex}] (aa) to node {} (aaa);
\draw[->, transform canvas={yshift=+0.5ex}] (aa) to node {} (aaa);

\draw[->] (b) to node {} (ba);

\draw[->,dashed] (a) to node {} (b);
\draw[->,dashed] (b) to node {} (c);

\draw[->, transform canvas={xshift=-0.25ex}, dashed] (aa) to node {} (ba);
\draw[->, transform canvas={xshift=+0.25ex}] (aa) to node {} (ba);

\draw[->] (aa) to node {} (b);
\draw[->] (aaa) to node {} (ba);
\draw[->] (ba) to node {} (c);

\draw[->, transform canvas={xshift=-0.2ex,yshift=-0.2ex}] (a) to node {} (ba);
\draw[->, transform canvas={xshift=+0.2ex,yshift=0.2ex}] (a) to node {} (ba);
\end{tikzpicture}
\] 
\caption{A sketch of the category $\D$, only objects $(a_0, \ldots, a_i)$ with $\sum a_i \leq 3$ drawn.
Marked arrows are dashed and arrows from or to $(1,2)$ omitted for readability (they are as for $(2,1)$).
We see that the top line is a copy of $\Delta_+$.} \label{fig:D}
\end{figure}
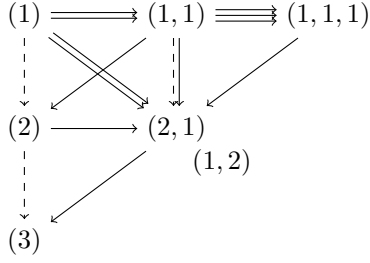

$\D$ has $\deltplus$ as a full subcategory consisting of all the lists $(1, \ldots, 1)$.
We can picture a part of $\D$ as shown in \cref{fig:D}. 

To every object of $\D$, we assign a degree by 
\begin{align}
 \deg(a_0,\ldots,a_m) &\defeq \left(2 \cdot \sum a_i\right) - (m+2).
\end{align}
Let us check that every non-identity morphism $f \in \D\left((a_0, \ldots, a_m),(b_0,\ldots,b_n)\right)$ increases the degree.
For this, let us write $I \subset [n]$ for the image of $f$ and $J$ for its complement, from which we get $\sum_{j\in I} b_j + \sum_{j \in J} b_j = \sum b_j$.
From the property of $f$, we get $\sum_{j \in I} b_j \geq \sum a_i$.
It is therefore sufficient to check the inequality $2 \cdot \sum_{j \in J} b_j > n - m$.
For $m > n$, this is immediate.
For $n > m$, we have that $f$ cannot be surjective, and we have $\sum_{j \in J} b_j \geq |J| \geq n - m \geq 1$, thus adding the factor $2$ on the left-hand side makes the inequality hold. 
For $n = m$, the inequality holds if $|J| \geq 1$, but if $|J| = 0$, then $f$ is injective and surjective, hence the identity.

We also see that there are only a finite number of objects of any degree, and every object is the codomain of a finite number of arrows in $\D$.
This means that the opposite $\oppo \D$ is an inverse category, which we call the \emph{inverse replacement} of $\oppo \Delta$.

\begin{definition}
 A \emph{simplicial type} is a Reedy fibrant diagram $S \co \oppo \D \to \type$ mapping each marked morphism in $\oppo \D$ to an equivalence.
\end{definition}

Let us describe our intuition for this construction.
Given a Reedy fibrant diagram $A \co \oppo \D \to \type$, the objects of the form $A_{(1,1,\ldots,1)}$ form a semisimplicial type, while the other objects encode the degeneracies.
For example, given a point $x : A_{(1)}$ and a loop $l : A_{(1,1)}(x,x)$ around this point, we have a type $A_{(2)}(a,l)$.
We can think of an element of this type as a proof that $l$ is the degenerated line we should get from $x$.
The marking ensures that the (only) projection $A_{(2)} \to A_{(1)}$ is an equivalence, which implies that there is exactly one such degenerated line $l$ for every point $x$.

A remark is that having such a type family $A_{(2)}$ with the property that the projection to $A_{(1)}$ is an equivalence is just as good as (i.e., equivalent to) having a function $\Pi (x : A_{(1)}), A_{(1,1)}(x,x)$. In \cref{fig:D}, we can see that this function can be constructed by inverting the marked arrow and composing with the single arrow $(2) \to (1,1)$.

Similarly, if we have two points $x,y : A_{(1)}$, a line $l : A_{(1,1)}(x,x)$, a proof $p$ that this line is the degeneracy of $x$, a second line $k : A_{(1,1)}(x,y)$, and a triangle filler $t : A_{(1,1,1)}(x,x,y,l,k,k)$, then $A_{(2,1)}(x,y,l,k,p,t)$ can be understood as the type of proofs that $t$ is the degeneracy of $k$ where $x$ is duplicated. Again, the relevant marking ensures that there is exactly one such degeneracy per line. The other degeneracy (where $y$ is duplicated) is induced by $A_{(1,2)}$.

If we have just a single point $x$ and we degenerate it to a line $l$, we have two ways to further degenerate it to a triangle (corresponding to the two degeneracy morphisms $[1] \to [2]$ in $\oppo \Delta$), but these two need to coincide.
This requirement is encoded by the projection $A_{(3)} \to A_{(1)}$ being an equivalence.

Similarly, any of the usual generating degeneracy maps $X_{[n-1]} \to X_{[n]}$ comes from morphisms in $\oppo \D$ of the form 
\begin{equation*}
 \underbrace{(1,\ldots,1)}_{n \text{ ones}} \dashleftarrow \underbrace{(1,\ldots,1,2,1,\ldots,1)}_{(n-1) \text{ ones, a single two somewhere}} \to \underbrace{(1,\ldots,1)}_{n+1 \text{ ones}}
\end{equation*}
where the left arrow is marked and can thus in the diagram be inverted.
The lists with higher numbers encode what is usually called the \emph{simplicial identities}, together with all coherence laws between them.

If we choose a number $n$ (externally) and consider the subcategory of $\Delta$ containing the objects $[0]$, $[1]$, \ldots, $[n]$, 
we can define finite versions of $\D$ and $\oppo \D$ consisting only of those lists with sum of the entries not exceeding $n+1$.
This can then be done inside type theory, i.e.\ we can formalise simplicial types up to level $n$ in proof assistants such as Agda, Coq, or Lean.

\subsection{Constructing Diagrams Over (Certain) Reedy Categories} \label{sec:ReedyDiags}

The construction of $\D$ is the special case of a more general construction for a Reedy category $R$. 
We denote maps in $\plus R$ by $x \cof y$ and maps in $\minus R$ by $x \fib y$.

\begin{definition}
 For a Reedy category $R$, 
 we define $D(R)$ to be the following category.
 Objects are arrows in $\minus R$, and a morphism between $s : x \fib y$ and $t : z \fib w$ is any morphism $f \in R(y,w)$ such that there exists a morphism $x \cof w$ in $\plus R$ which makes the following square commute:
\[
\begin{tikzpicture}[x=1.7cm,y=-1.2cm]
 \node(x) at (0,0) {$x$};
 \node(y) at (0,1) {$y$};
 \node(z) at (1,0) {$z$};
 \node(w) at (1,1) {$w$\rlap{.}};
 
 \draw[->>,left] (x) to node {s} (y);
 \draw[->>,right] (z) to node {t} (w);
 \draw[>->,dotted] (x) to node {} (z);
 \draw[->,above] (y) to node {f} (w);
\end{tikzpicture}
\]
Note that the dotted arrow is only required to exist, it is not part of the structure of the morphism $f$.
We say that a morphism in $D(R)$ is \emph{marked} if $f$ is an identity in $R$.
\end{definition}

For the special case of $\Delta$, the maps in $\minus \Delta$ are the surjections, and a surjection $s : [m] \fib [n]$ can be represented by as a list $(a_0,\ldots,a_n)$, with $a_i$ being the size of the preimage of $i$ under $s$.
Under this translation, the existence of a lift in $\deltplus$ corresponds exactly to the condition we had put on morphisms.

\begin{remark}
Our construction of $\D$ from $\Delta$ is closely related to the ``fat $\Delta$'' of \cite{kock-weak-identity-arrows}, which has the same objects as our $\D$, but morphisms for which the required lift of a map in $\Delta$ to $\deltplus$ mentioned above is part of the data.
Kock uses it to define so-called fair categories, a tool for attacking Simpson's conjecture that composition and exchange laws can be made strict in higher categories, leaving only weak unit laws.
\end{remark}

We can see that the arrows in $D(R)$ go only in one direction as follows.
Assume we have an infinite chain:
\[
\begin{tikzpicture}[x=-1.2cm,y=-0.8cm]
 \node(x0) at (0,0) {$x_0$};
 \node(y0) at (0,1) {$y_0$\rlap{.}};
 \node(x1) at (1,0) {$x_1$};
 \node(y1) at (1,1) {$y_1$};
 \node(x2) at (2,0) {$x_2$};
 \node(y2) at (2,1) {$y_2$};
 \node(x3) at (3,0) {\ldots};
 \node(y3) at (3,1) {\ldots};
  
 \draw[->>] (x0) to node {} (y0);
 \draw[->>] (x1) to node {} (y1);
 \draw[->>] (x2) to node {} (y2);

 \draw[->] (y1) to node {} (y0);
 \draw[->] (y2) to node {} (y1);
 \draw[->] (y3) to node {} (y2);
 
 \draw[>->] (x1) to node {} (x0);
 \draw[>->] (x2) to node {} (x1); 
 \draw[>->] (x3) to node {} (x2);
\end{tikzpicture}
\]

Because the upper horizontal maps can only increase the degree in $R$, they are eventually (say, for all indices greater or equal to $N$) all identities.
Thus, for $k > N$, we have that the composition $x_{k+1} \fib y_{k+1} \to y_k$ is in $\minus R$.
If we factor the map $y_{k+1} \to y_k$ as $y_{k+1} \fib \hat y_k \cof y_k$, we see from the uniqueness of the factorisation of $x_{k+1} \to y_k$ that the map $\hat y_k \cof y_k$ has to be the identity. In summary, the lower horizontal maps after index $N$ are all in $\minus R$ and decrease the degree in $R$, which however is bounded by $\deg(x_N)$, which means that eventually, \emph{all} horizontal morphisms become identities.

For $R \jdeq \Delta$, we can 
define a degree function on $D(R)$ which uses the degree function on $R$ by assigning $s : x \fib y$ the value $2\cdot \deg(x) - \deg(y)$.
For general $R$, this does not necessarily work, but an appropriate degree function can still be defined.

Just as in the special case, we can now consider Reedy fibrant diagrams over $\oppo{D(R)}$ with the condition that marked arrows get mapped to equivalences.
If the canonical functor $D(R) \to R$ is an opfibration, 
then this gives a well-behaved notion of diagrams equivalent to homotopy coherent diagrams with identities.
This is in particular the case if $R$ is $\Delta$, or any ``finite version'' of $\Delta$ (i.e.\ $\Delta$ restricted to objects $[0], \ldots, [n]$).
The details will be given in the full version of the paper, but some core ideas are explained in the next section.

\section{Equivalences between the constructions} \label{sec:equivalencesBetweenConstructions}

In~\cite{we_futurework}, we will present full proofs of more general versions of the following results: 
\begin{theorem}
 The notions of different diagrams that we have presented in this short paper are equivalent whenever it makes sense.
 In detail:
 \begin{enumerate}
  \item For an inverse category $I$, the strictification construction of \cref{sec:strictification} establishes an equivalence between weak and strict diagrams over $I$ as long as both types exist (e.g.\ if $I$ is finite). 
  \item The notion of a homotopy coherent diagram over $I$ and a general homotopy coherent diagram with identities over $I$ are equivalent, in the sense that the projection (fibration) mentioned in \cref{sec:cohId} has contractible fibres. (This needs that the type theory supports the occurring infinite notions.)
  Both notions thus coincide with that of Reedy fibrant diagrams.
  \item For Reedy categories satisfying certain technical conditions (which ensure that $D(R) \to R$ is an opfibration, and which hold e.g.\ for $\Delta$), general homotopy coherent diagrams with identities are equivalent to the diagrams over Reedy categories as constructed in \cref{sec:ReedyDiags}. 
 \end{enumerate}
\end{theorem}
We further hope to be able to show even stronger connections, as follows. 
Using the notion of a higher category or complete semi-Segal type as suggested in~\cite{altenCapKraus_infCats}, we can consider a complete semi-Segal type (or an $\infty$-semicategory) of diagrams in each case. 
The equivalences we establish should then not only be equivalences between types of diagrams, but rather higher equivalence between $\infty$-semicategories.
In the remaining part of this short paper, we give \emph{very} brief sketches of the arguments to (hopefully) make some key ideas understandable.

\subsection{General techniques: the Segal condition and inner, left, and right fibrations} \label{sec:someFinalGeneralTechniques}

Let $A \co \deltop \to \type$ be a semisimplicial type.
Recall that we think of $A_{[n]}$ as an $n$-dimensional tetrahedron, with $2^n-1$ many cells (starting with $n+1$ points, $\binom{n+1}{2}$ lines, and so on).
We can consider the usual constructions of simplicial sets: for example, if we remove the single $n$-dimensional cell from $A_{[n]}$ and one of the $(n-1)$-dimensional cells, we get what is called a \emph{horn}, for which we could write $A(\Lambda^{[3]}_i)$, where $i$ is the number of the vertex opposite to the removed $(n-1)$-dimensional cell.
Some semisimplicial types $A$ have the property that projections of the form $A_{[n]} \fib A(\Lambda^{[3]}_i)$ are equivalences, and in this case, we say that $A$ has \emph{contractible horn fillers} or that it is \emph{Kan fibrant}.
In the case of $\ty$, we observe that we have contractible fillers for all inner horns (i.e.\ $0 < i < n$).
In fact, having contractible fillers for inner horns is equivalent to having the type-theoretic \emph{Segal condition}, where the projection is an equivalence that maps the full tetrahedron to the sub-tetrahedron consisting only of the points and a chain of lines.
This holds for $\ty$ by construction.

The notion of a Reedy fibrant diagram is only the special case of a Reedy fibration with the terminal diagram as the codomain, and for Reedy fibrations, the notions of inner, Kan, left, or right fibrations all make sense and express that different selections of horns have contractible fillers. 
In particular, we can construct a Reedy fibration $\ty^\bullet \to \ty$, where $\ty^\bullet$ can be thought of a as a \emph{pointed} version of $\ty$.
This Reedy fibration serves as a \emph{left fibrations classifier} and can be constructed as a restricted version of Shulman's universe fibration $\tilde U \fib U$.
Another possible definition is to set $\ty^\bullet_{[n]}(\ldots) \defeq \ty_{[n+1]}(\unit, \ldots)$.
If we look at the strictification construction in \cref{sec:strictification}, we see that this \emph{index shifting} appears already there.

\subsection{Weak versus strict diagrams over inverse categories} \label{subsec:weakStrictEquiv}

The following is a proof sketch which gives intuition for the equivalence of Reedy fibrant and homotopy coherent diagrams over a finite inverse category $I$.
We found that this proof sketch is not the easiest to complete as some details are somewhat tricky; but the more elegant proof that will be presented in~\cite{we_futurework} is not suitable for providing intuition using only limited space.

We do induction on the index category $I$.
Thus, we assume that $i$ is some maximal object of $I$ (receives no non-identity arrows), and write $I_0$ for the category $I$ with $i$ removed.
By induction, the equivalence holds for the types of diagrams over $I_0$.

Let a homotopy coherent diagram $X$ over $I_0$ be given. 
To extend this to a homotopy coherent diagram over $I$, we need one new component (context entry) of the form $\left(X_{i \to \ldots} : \ty_{[k]}(\ldots)\right)$ for every sequence of composable morphisms starting with $i$, with $k$ being the length of this sequence.
We can split this and say that the new components are a single type $(X_i : \UU)$ and one component $\left(X_{i \to \ldots} : \ty_{[k]}(X_i, \ldots)\right)$ for every non-zero sequence starting with $i$ (note that the first component of the matching object will always be $X_i$).
It is intuitive (caveat: but not easy to show) that $\ty_{[k]}(X_i, \ldots)$ is equivalent to the function type $X_i \to \ty_{[k]}(\unit, \ldots)$, where the second occurrence of $\ldots$ already makes use of this equivalence on lower levels (it may be helpful to spell this out explicitly for a few low values of $k$).
If we now look at the strictification construction in \cref{sec:strictification}, we see that the matching object $M^X_i$ we construct there consists exactly of these components in $\ty_{[k]}(\unit, \ldots)$.
In short, the new components needed to extend the diagram over $I_0$ to a diagram over $I$ can be represented as a type $X_i$ together with a function $X_i \to M^X_i$. 
It is a general principle that, given a type $M$, the type $\Sigma (A : \UU). (A \to M)$ is equivalent to $M \to \UU$. 
What this means here is that the new components are equivalent to a type family $M^X_i \to \UU$, and this is exactly what is needed to extend the Reedy fibrant diagram, constructed from $X$ over $I_0$, to a Reedy fibrant diagram over $I$.

A possibly helpful exercise is to use this strategy explicitly to prove that the two types of diagrams over $\E$, given in \cref{ex:E-Reedyfibrant-internally,ex:weak-diagram}, are equivalent.

\subsection{Neutrality of identities} \label{subsec:neutralityOfIds}

For an inverse category, we claim that the two versions of homotopy coherent diagrams over $I$ (\cref{def:hcoh,def:hcohWithId}) are equivalent.
Already the special case that the inverse category is the terminal category $\term$ with a single object and no non-identity morphism is interesting.
$\int \posnerv \term$ is still $\term$, and a homotopy coherent diagram is thus just given by a single type. 
However, $\int \fullnerv \term$ is infinite, and a general homotopy coherent diagram $h$ has one component for each sequence of identity arrows, i.e.\ one component $h_n$ for each $n \in \N$, starting as follows. 
$h_0$ is just a type.
$h_1 : h_0 \to h_0$ is a function, which our \emph{Harpaz condition} requires to be an equivalence.
Further, we have $h_2 : h_1 \circ h_1 = h_1$, then $h_3$ stating that the different ways of composing the equality $h_2$ with itself coincide, then $h_4$ certifying the coherence of $h_3$, and so on.

Let us write $H_n$ for the limit over the subcategory of $\int \fullnerv \term$ consisting of only those objects of degree $\leq n$;
this means, $H_n$ will consist of only finitely many components $h_0$, \ldots, $h_n$.
The canonical projection $H_\infty \fib H_0$ is the map described in \cref{sec:cohId}, i.e.\ the map from weak diagrams with identities to those without identities, which we claim to be an equivalence.
To show this, we observe that $H_\infty$ is the limit of
\begin{equation}
 H_0 \lfib H_2 \lfib H_4 \lfib H_6 \lfib \ldots,
\end{equation}
and we show that every map in this sequence is an equivalence.
To show that $H_{2n+2} \fib H_{2n}$ is an equivalence, we show that the components $h_{2n+1}$ and $h_{2n+2}$ together form a contractible pair.
In the case $n=0$, this is an auto-equivalence $h_1 : h_0 \simeq h_0$ and a proof of $h_1 \circ h_1 = h_1$, and this pair is easily seen to be contractible.
As it was pointed out to us by Peter LeFanu Lumsdaine, this is a type-theoretic version of what is known as \emph{dunce's hat}~\cite{zeeman1963dunce}: \emph{eunce's hat} in topology is the simplest example of a space which is contractible but not collapsible.
In type-theoretic terms, this means we have a type which is contractible but not a ``singleton'' or a collection of ``singletons'' (types of the form $\Sigma (a : A). a = a_0$, which are known to be contractible).

If the two component $h_{2n+1}$ and $h_{2n+2}$ formed a ``horn filler'' in $\ty$, they would be contractible (see \cref{sec:someFinalGeneralTechniques}), but this is not the case as the type of $h_{2n+2}$ has $(2n+3)$ occurrences of $h_{2n+1}$ inside instead of a single one.
However, as we will show in our future article, the odd number of occurrences ensures that the pair is contractible nevertheless.
The general case of an inverse category $I$ instead of just $\term$ is only slightly more difficult.

\subsection{Simplicial types versus general homotopy coherent diagrams over the simplex category} \label{subsec:simplTypesVsHC}

Simplicial types are by definition Reedy fibrant diagrams over $\oppo \D$ respecting the markings.
The previous result 
sketched in \cref{subsec:weakStrictEquiv} 
can be extended to this marked case, and we can thus equivalently consider general homotopy coherent diagrams over $\oppo \D$ with some conditions.
There is a canonical functor $F \co \oppo \D \to \oppo \Delta$, which turns out to be a Grothendieck fibration. 
This allows us to find, given any sequence of morphisms in $\oppo \Delta$ (i.e.\ any object of $\int \fullnerv \oppo \Delta$),
a sequence in $\oppo \D$ lying over it. 
Using techniques described in \cref{sec:someFinalGeneralTechniques,subsec:neutralityOfIds}, we then show that all the additional (``non-minimal'') sequences in $\oppo \D$ cancel each other out when considering general homotopy coherent diagrams.

\section*{Acknowledgments}
We would like to thank Paolo Capriotti for many discussions and for pointing out Harpaz' trick, without which the presented notion of diagrams in \cref{sec:cohId} would be significantly more involved or nonexistent.
We are also grateful to Thorsten Altenkirch, Steve Awodey, Ulrik Buchholtz, and Peter LeFanu Lumsdaine for their comments on this work.

\bibliographystyle{plain}
\bibliography{master}

\end{document}